# *The Numerical Solution in the Sense of Prager&Synge*

A.K. Alekseev[1,2], A.E. Bondarev[3]

[1]Moscow Institute for Physics and Technology, Moscow, Russia
[2]RSC Energia, Korolev, Russia
[3]Keldysh Institute of Applied Mathematics RAS, Moscow, Russia

e-mail: alekseev.ak@phystech.edu, bond@keldysh.ru

***Abstract***

The solution in sense of Prager&Synge is the alternative to the commonly used notion of the numerical solution, which is considered as a limit of grid functions at mesh refinement. Prager&Synge solution is defined as a hypersphere containing the projection of the true solution of the system of partial differentiation equations (PDE) onto the computational grid and does not use any asymptotics. In the original variant it is determined using orthogonal properties specific for certain equations. In the proposed variant, the center and radius of the hypersphere is estimated using the ensemble of numerical solutions obtained by independent algorithms. This approach may be easily expanded for solutions of an arbitrary system of partial differentiation equations that significantly expands the domain of its applicability.

Several options for the computation of the Prager&Synge solution are considered and compared herein. The first one is based on the search for the orthogonal truncation errors and their transformation. The second is based on the orthogonalization of approximation errors obtained using the defect correction method and applies a superposition of numerical solutions. These options are intrusive.

In third option (nonintrusive) the information regarding orthogonality of errors, which is crucial for the Prager&Synge approach method, is replaced by information that stems from the properties of the ensemble of numerical solutions, obtained by independent numerical algorithms. The values of the angle between the truncation errors on such ensemble or the distances between elements of the ensemble may be used to replace the orthogonality. The variant based on the width of the ensemble of independent numerical solutions does not require any additional *a priori* information and is the approximate nonintrusive version of the method based on the orthogonalization of approximation errors.

The numerical tests for two dimensional inviscid flows are presented that demonstrate the acceptable effectivity of the approximation error estimates based on the solution in the Prager&Synge sense.

## 1. Introduction

The verification of numerical solutions is of the current interest due to the great practical demand for the reliable computations that is stated by the modern standards [1,2]. Herein, we discuss the approximation error estimation that is the necessary component of the verification using a notion of the Prager&Synge numerical solution.

We denote an operator form of the system of equations by $A\tilde{u} = \rho$ . The discrete operator, which approximates this system on some grid, is noted by $A_h u_h = \rho_h$. The Lagrange identity $(A_h u_h, v)_{L_2} = (u_h, A_h^* v)_{L_2}$ defines the adjoint operator $A_h^*$. The numerical solution $u_h$ is a grid function (vector $u_h \in R^M$ , $M$ is the number of grid points), $\tilde{u}_h \in R^M$ is the projection of a true solution onto the considered grid ($h$ is the step of discretization), $\Delta\tilde{u}_h = u_h - \tilde{u}_h$ is the approximation error, $\Delta u_h$ is an estimate of this error.

The *a priori* error estimates, widely used at design and analysis of numerical schemes convergence, have the form $\|\Delta\tilde{u}_h\| \le Ch^n$ ($n$ is the order of approximation). These estimates have the general form and do not depend on the specific solution. Unfortunately, the unknown constant $C$ hinders the practical applications of this approach.

The *a posteriori* error estimate has the form $\|\Delta\tilde{u}_h\| \le I(u_h)$ containing some computable error indicator $I(u_h)$, which depends on the numerical solution $u_h$ (that restricts the generality of results by the concrete solution), but has no unknown constants. So it provides a computable measure of the discretization error for the analyzed numerical solution and is highly interesting from the practices standpoint.

The *a posteriori* error analysis has the long history (starting from [3,4,5,6]). At present, the high efficiency of this technique [7,8,9] is achieved in the domain of the finite-element analysis (due to smooth enough solutions of elliptic and parabolic equations). In contrast, the progress in the *a posteriori* error estimation is limited for the problems governed by the equations of the hyperbolic or mixed type (typical CFD problems). This situation is caused by the discontinuities (shock waves, contact surfaces) that are specific for these problems. The modern standards [1,2] recommend the Richardson extrapolation (RE) [11,12] for the error estimation in CFD applications. However, RE is the essentially semiheuristic approach, since it is based on the leading error term asymptotics. For the compressible flows, it is restricted by the variation of the convergence order over the flowfield [13, 14]. As the cure, the generalized Richardson extrapolation (GRE) [15,16] may be applied that provides an estimate for the distribution of the convergence order over grid. Unfortunately, GRE demonstrates a very high level of oscillations in estimates of the convergence order and entails the great computational burden (it requires at least four consequent refinement of mesh [15,16]).

So, the computationally inexpensive estimators of the approximation error in CFD domain are of current interest and we consider, herein, some new options related with the famous Prager&Synge [3,4,5] method. From the historical perspective, the technique for the *a posteriori* error has started exactly from this method. At present, *a posteriori* error estimation is based on other ideas [6,7,8]. Nevertheless, the Prager&Synge method seems to be highly underestimated both from the viewpoint of general idea and from the viewpoint of applicability domain.

The original notion of the numerical solution is the key idea of Prager@Synge method. In contrast to the standard numerical solution realized by a sequence converging to the true solution at grid step diminishing, the solution in Prager@Synge sense is some hypersphere that contains the true solution. This approach is an ideal one from the viewpoint of verification. It is interesting that it is not based on a mesh refinement asymptotics.

Several options are considered and compared for the estimation of the numerical solution in Prager@Synge sense. All of them apply an ensemble of numerical solutions, obtained by independent algorithms. These algorithms may differ in the order of convergence or in the structure (for example, based on the Riemann problem solution or on the finite differences). We consider three main variants for the estimation of the Prager&Synge solution.

First, we construct an artificial truncation error orthogonal to the estimate of the truncation errors performed by some numerical algorithm using special adjoint postprocessor.

Second, we search for the superposition of numerical solutions (obtained by independent algorithms) that provides the approximation error, orthogonal to the certain solution error.

Third, we demonstrate that several methods for the approximation error norm estimation, whish are based on the triangle inequality [17,18,19], the width of the ensemble of solutions [20], and the angle between approximation errors [21], may be considered as some semiheuristical nonintrusive approximations of the Prager@Synge method.

The two dimensional compressible flows, governed by the Euler equations, are computed in the numerical tests in order to verify and compare the above mentioned algorithms. These flow patterns contain shocks and contact discontinuities that provides most rigorous environment for the approximation error estimation. The analytical solutions for Edney-I and Edney-VI flow patterns [22] are used to obtain the "true" discretization error by comparing with the numerical solutions.

Despite the fact that our numerical tests concern CFD domain, we believe that all methods are applicable for the numerical solutions for a wide range of the PDE systems.

The paper is organized as follows. In Section 2 we remind the classical variant of the Prager&Synge method and discuss the prospects for the expansion of its applicability domain. The Section 3 discusses the universal numerical solution in the Prager&Synge sense. Section 4 surveys the methods proposed for a computation of the solutions in the Prager&Synge sense. Section 5 describes the test problems and the set of the numerical algorithms used to generate the ensemble of solutions. The results of the Prager&Synge solution estimation by different approaches are presented in Section 6. The discussion of results is provided in Section 7. The conclusions are presented in Section 8.

## 2. Original Prager&Synge method

The "Prager&Synge method" developed in [3,4,5] and noted as "hypercircle method", is historically the first example of a method for *a posteriori* error estimation. This method was the result of the collaboration of German specialist in applied mathematics (William Prager) and Irish mathematician and physicist (well known as the specialist in the interior of the black holes) John Lighton Synge, who supported and promoted this method for a quarter of a century. Perhaps, this unique collaboration is the reason for highly unusual nature and properties of the "hypercircle method". The term "hypercircle" is currently used in distant science domains, so, in order to avoid confusion, we mark this method by the "Prager&Synge method". At present, the Prager&Synge method is applied in rather narrow domain of applications. By this reason, modern publications [23-26] give a limited presentation of deep ideas, described by old works. So, we provide herein some brief presentation of Prager&Synge method with illustrations and citations.

The Prager&Synge method is commonly applied for the solution of the Poisson equation

$$\nabla^2 u = \rho\,,\ u_\Gamma = f\,. \tag{1}$$

Two additional linear subspaces of functions are used:

the subspace of gradients $\partial v / \partial x_i$ of functions $v \in H^1$ that satisfy boundary conditions

$$v_\Gamma = f\,, \tag{2}$$

the subspace of functions $q \in H^1$, such that

$$div(q) = \rho\,. \tag{3}$$

The following orthogonality condition holds for these functions and exact solutions of the Poisson equation $\tilde{u}$

$$((\nabla v - \nabla \tilde{u}),(q - \nabla \tilde{u})) = 0\,. \tag{4}$$

Usually $v$ is chosen as the numerical (finite elements) solution $v = u_h$. All relations are valid both for the function space and for the finite-dimensional case of the grid functions. All norms and scalar products correspond to $L_2$ and its finite-dimensional analogue.

The main idea of this approach is purely geometrical. We provide the vast citation from [5], since we are unable to compete with authors of the method in vividness and lucidity.

*"Imagine a flat piece of country. Imagine in it two strait narrow roads which intersect at right angles. Imagine two blind men set down by helicopter, one on each road. Provide them with a long measuring tape, stretching from the one man to the other man. Let the tape be embossed, so they can read it, even though blind. Without moving, let them tighten the tape and read it. Suppose it reads 100 yards apart.*

***We are to ask each of them how far he is from the cross roads.***

*These men, though blind, are very intelligent, and they know, that the hypotenuse of a right-angled triangle is greater than either the other two sides. So each of the men says that his distance from the cross-roads is less than (or possible equal to) 100 yards".*

Fig. 1 ([5]) illustrates the "Two blind men" problem.

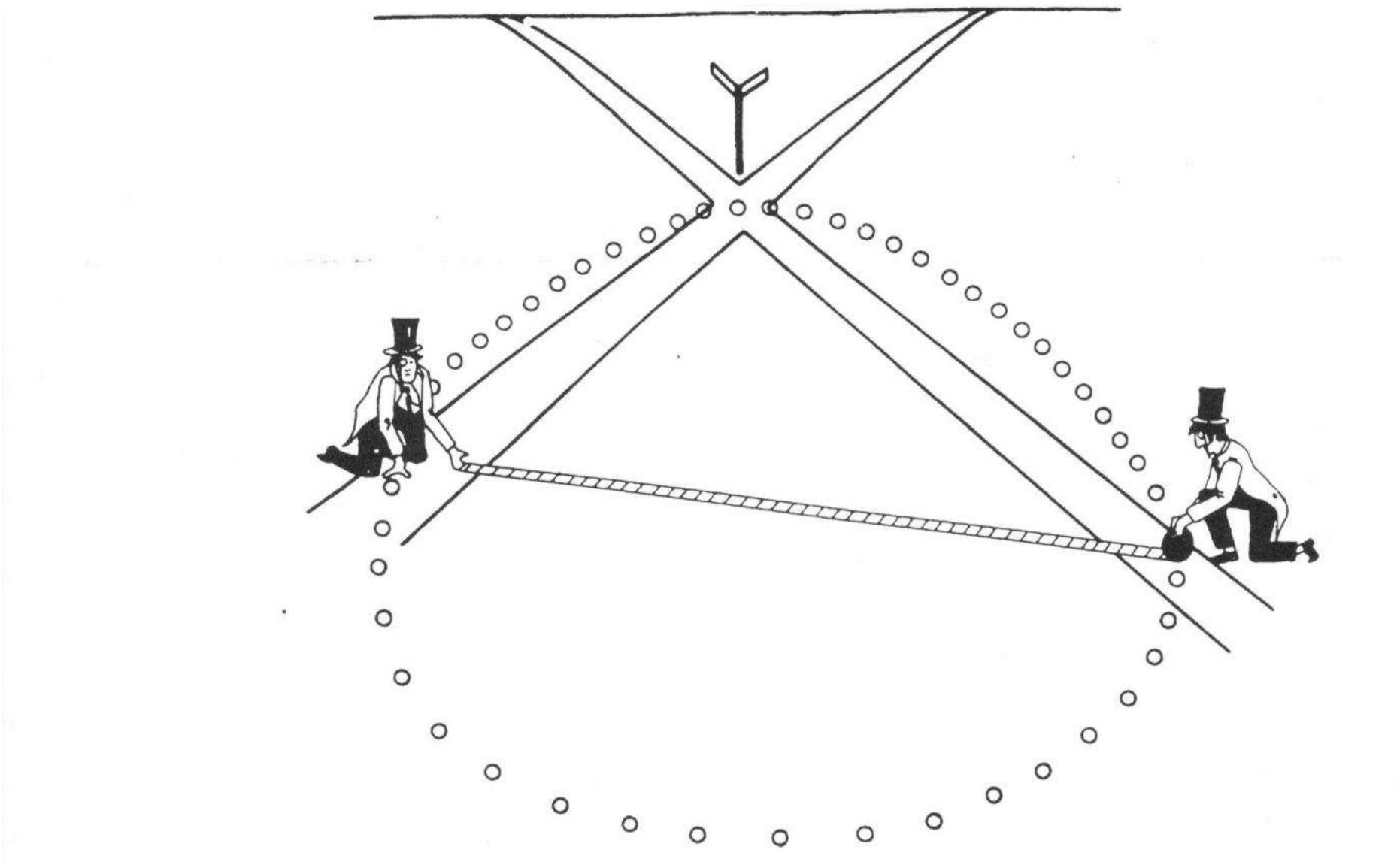

Fig. 1. Two blind men and their hypercircle from [5]

Thus, the Prage@Synge method provides *a posteriori* error estimate from ***purely geometrical ideas*** and without any unknown constants that cause the great interest to it both from the theoretical and practical viewpoints.

The orthogonality relation (Eq. (4)) engenders the inequality for the approximation error (distance between exact and numerical solutions), since the leg of a right triangle is less the hypotenuse

$$\|\nabla(\tilde{u} - u_h)\| \le \|\nabla u_h - q\| . \tag{5}$$

It should be mentioned that the Prager&Synge method has significant restrictions:

1. The domain of its applicability is limited by the problems governed by the Poisson equation and some closely related problems (biharmonic equation [23]). In general case the orthogonality relation may be stated as $(Bv - B\tilde{u})(q - B\tilde{u}) = 0$. It may be recast as

$$(Bv - B\tilde{u})(q - B\tilde{u}) = (v - \tilde{u})(B^* q - B^* B\tilde{u}) = 0 . \tag{6}$$

It holds for equations $A\tilde{u} = \rho$ at $A = B^* B$, $\rho = B^* q$, $q = B^{*-1}\rho$. This provides a bit more general form of Eq. (5)

$$\|B(\tilde{u} - u_h)\| < \|Bu_h - q\| . \tag{7}$$

So, the Prager&Synge method may be applied to very narrow domain of equations, which does not include most of practically interesting systems of equations (such as Euler on Navier-Stokes).

2. The estimation of the auxiliary function $q$, defined by the equation (3), entails numerical difficulties.

3. By this method one estimates the norm of the error of the solution gradient (or expression $B\tilde{u}$) and is unable to estimate the error norm $\|\tilde{u}-u_h\|$, which is of the primal interest from the practitioners needs**.**

Perhaps by these reasons, the mainstream of *a posteriori* error estimates is, at present, based on the quite different ideas [6,7]. Nevertheless, the approach by Prager&Synge is highly underestimated [27] and may be significantly expanded that we try to demonstrate below.

## 3. The numerical solution in the sense of Prager&Synge

The main idea by Prager&Synge consists in the original notion of the numerical solution. We mark it as the numerical solution in the Prager&Synge sense. The numerical solution, defined by the standard way, is considered to be an element of the sequence $u_{h_m}$ that is assumed to converge to the projection of the exact solution $u_{h_m} \to \tilde{u}_{h_m}$ as the step of discretization decreases ($h_m \to 0$ at $m \to \infty$). The fundamental theorem of numerical analysis (Lax-Richtmyer Equivalence Theorem [10]) defines the conditions at which such sequence converges to the exact solution. From the standpoint of practice the elaborated processes of the mesh adaptation and refinement are the key elements enabling the success of this approach.

Contrary to the standard approach, Synge stated ([4], p. 97): *"In general, a limiting process is not used, and we do not actually find the solution.... But although we do not find it, we learn something about its position, namely, that it is located on a certain hypercircle in function space"*.

In finite dimensional case the "hypercircle" by Prager&Synge corresponds to the hypersphere. So, the numerical solution by Prager&Synge is considered as a hypersphere with the center $C_h = u_h$ and a radius $R_h$ containing the projection of the true solution onto the grid $\|\tilde{u}_h - C_h\| \le R_h$. The concept of the numerical solution by Prager&Synge is closely related to the purposes of the *a posteriori* error estimation. In the contrast to the standard approach, the convergence of numerical solution is not obligatory ***("a limiting process is not used")***. The mesh refinement is not mandatory also and should be performed only, if the magnitude of $R_h$ is not acceptable from the viewpoint of practical needs (that may be checked using Cauchy–Bunyakovsky–Schwarz inequality for valuable functionals (some examples are presented by [20])). In this way, the natural stopping criterion for the mesh refinement termination may be stated.

So, the solution in sense of Prager&Synge enables to avoid the "tyranny" of the mesh refinement and adaptation that governs the modern CFD applications.

As we already discussed, the standard domain of the Prager&Synge method applicability is rather narrow and does not include the Euler or Navier-Stokes equations that are of the main interest in CFD. Nevertheless, we believe that the domain of the Prager&Synge method applicability may be extended to the arbitrary PDE system. In the spirit of the Prager&Synge solution we search for some auxiliary solution $u_\perp$ such that

$$(u_h^{(i)} - \tilde{u}_h, u_\perp - \tilde{u}_h)_\Omega = (\Delta\tilde{u}_h^{(i)}, \Delta\tilde{u}_\perp)_\Omega = 0\,. \qquad (8)$$

First, we consider the application of this approach for the set of numerical solutions (grid functions) obtained by independent numerical methods on the same grid. Let $u_h^{(1)} \in R^M$ and $u_h^{(2)} \in R^M$ be the numerical solutions having the orthogonal approximation errors $\Delta u_h^{(i)} = u_h^{(i)} - \tilde{u}_h$ ($\Delta u_h^{(1)} \perp \Delta u_h^{(2)}$). The Prager@Synge method, corresponding to the ideal case of the orthogonal approximation errors (the angle between errors $\alpha = \pi/2$), is illustrated by Fig. 2.

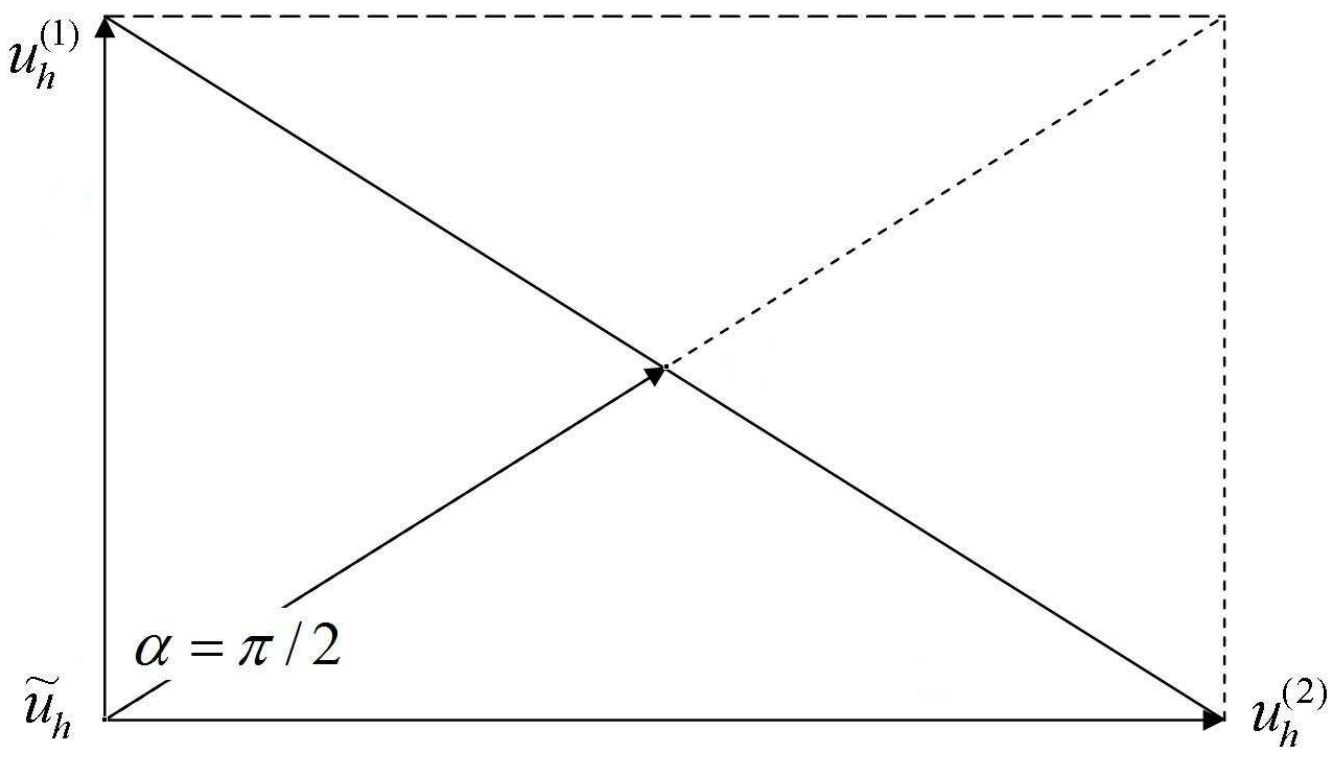

Fig. 2. $\alpha = \pi/2$

The inequality

$$\left\|\tilde{u}_h - u_h^{(k)}\right\| < \left\|u_h^{(1)} - u_h^{(2)}\right\|, k = 1,2 \tag{9}$$

is valid, since any leg of the right triangle is less the hypotenuse. Herein, we assume $\Delta u_h^{(2)} = \Delta u_\perp$ for the illustration, that is not correct in general, since approximation errors $\Delta u_h^{(i)}$ are not orthogonal as demonstrated by [21].

## 4. Methods for the approximation of numerical solution in Prager&Synge sense

Herein we consider some algorithms for the search for the Synge numerical solutions. All of them are based on the usage of the ensemble of numerical solutions obtained by independent algorithms.

### 4.1. The orthogonal truncation error based method

Truncation errors $\delta u_h^{(i)}$ may be estimated either by the differential approximation [28] or by the special postprocessor [29]. In accordance with [20,21] the angles between truncation errors $\delta u_h^{(i)}$ on the ensemble of solutions are close to $\pi/2$ (orthogonal). Certainly, the estimates of truncation error may be directly used for calculation of the approximation error $\Delta u_h^{(i)}$

$$\Delta\tilde{u}_h^{(i)} \approx \Delta u_h^{(i)} = A_h^{-1}\delta u_h^{(i)} \tag{10}$$

by solving a special problem for disturbances (defect correction problem, for instance [30]).

But it is more interesting to restate the orthogonality condition (Eq. (8)) in terms of nonintrusively computable truncation errors. Let's select some vector $\theta \perp \delta u_h^{(i)}$ and transform it using forward and adjoint operators:

$$\delta u_\perp = A_h A_h^* \theta\,. \tag{11}$$

Expression (11) may be inverted to obtain

$$\theta = A_h^{-1*} A_h^{-1} \delta u_\perp . \tag{12}$$

The orthogonality condition $\theta \perp \delta u_h^{(i)}$ is equivalent to relation (8) since

$$(\delta u_h^{(i)}, \theta)_\Omega = (\delta u_h^{(i)}, A_h^{-1*} A_h^{-1} \delta u_\perp)_\Omega = (A_h^{-1} \delta u^{(i)}, A_h^{-1} \delta u_\perp)_\Omega = (\Delta u_h^{(i)}, \Delta u_\perp)_\Omega = 0 , \tag{13}$$

if

$$\Delta u_\perp = A_h^{-1} \delta u_\perp = A_h^* \theta . \tag{14}$$

So, $\Delta u_\perp$ may be computed nonintrusively using postprocessor that implements the solution of the adjoint problem (similarly to [27,29]). The adjoint postprocessor may be constructed using any solver designed for solution of adjoint equations.

The design of $\theta$, which is orthogonal to $\delta u_h^{(i)}$, may be conducted as follows:

1. By choice of arbitrary $\theta \perp \delta u_h^{(i)}$ (under condition $\|\theta\| = \|\delta u_h^{(i)}\|$).
2. By choice of $\theta$ that is equal to the truncation error of certain additional numerical solution

$$\theta = \delta u_h^{(k)} . \tag{15}$$

The latter option is caused by the numerical tests [21] that demonstrate the truncation errors, corresponding to independent algorithms, to be close to orthogonal. The orthogonality can be numerically checked.

So, the following inequality may be obtained:

$$\|\Delta u_\perp\| \le \|A_h^*\| \cdot \|\delta u_h^{(k)}\| \tag{16}$$

Thus $\|\Delta u_\perp\| \to 0$ at $\|\delta u_h^{(k)}\| \to 0$ that demonstrates the convergence of the orthogonal disturbance at mesh refinement (at bounded $\|A_h^*\|$).

Having $\delta u_\perp = A_h A_h^* \theta$, we may compute the auxiliary solution

$$u_\perp = A^{-1}(\rho + \delta u_\perp) . \tag{17}$$

In result the following inequality holds

$$\|u_h^{(k)} - \tilde{u}\| \le \|u_h^{(k)} - u_\perp\| . \tag{18}$$

The present approach constructs some hypersphere containing the true solution. So, the centre $C_h = u_h^{(k)}$ and radius $R = \|u_h^{(k)} - u_\perp\|$ determine the solution in the sense of the Prager&Synge.

Unfortunately, this approach is extremely complicated form the algorithmic viewpoint (mainly due to the application of the adjoint solver) and unstable (due to the differentiation of the non regular function by the adjoint solver). The numerical tests demonstrated highly pessimistic estimations (too great error norm) for above considered method.

So, the above discussion illustrates only the feasibility of $u_\perp$ construction and has mainly the heuristic value.

### 4.2. The orthogonal approximation error based method

As the second option for the generation of the Prager&Synge solution, we construct the auxiliary solution $u_\perp$ using an orthogonalization of the approximation errors.

The solution $u_\perp$ is defined as the superposition of $N$ numerical solutions (obtained by independent algorithms on the same grid)

$$u_\perp = \sum_{i=1}^{N} w_i u_h^{(i)} = w_i u_h^{(i)}, \tag{19}$$

$$\sum_{i=1}^{N} w_i = 1. \tag{20}$$

First, we select some basic numerical solution $u_h^{(0)}$ (we use it as the centre of the hypersphere) and estimate corresponding approximation error $\Delta u_h^{(0)} = A_h^{-1}\delta u_h^{(0)} \approx \Delta\tilde{u}_h^{(0)} = u_h^{(0)} - \tilde{u}_h$ using defect correction approach. Second, we search for the auxiliary solution $u_\perp$ (Eqs. (19),(20)) that is specific by the approximation error $\Delta u_\perp$ orthogonal to the error of the basic numerical solution

$$(\Delta u_h^{(0)}, \Delta u_\perp) = 0. \tag{21}$$

We assume that some $\Delta u_\perp$ exists due to Eq. (14). The auxiliary solution $u_\perp$ defines the radius of the hypersphere

$$R = \left\| u_h^{(0)} - u_\perp \right\| \tag{22}$$

that contains a true solution. The centre $C_h = u_h^{(0)}$ and radius $R = \left\| u_h^{(0)} - u_\perp \right\|$ determine the numerical solution in the sense of the Prager&Synge.

Under conditions (Eqs. (19),(20)) the “approximation error” of the auxiliary solution is the superposition of true (unknown) errors of the ensemble of solutions $\Delta u_\perp = w_i \Delta\tilde{u}_h^{(i)}$.

We search for $\{w_i\}$ in the variational statement

$$\{w_i\} = \arg\min \varepsilon(\vec{w}), \tag{23}$$

$$\varepsilon(\vec{w}) = (\Delta u_h^{(0)}, w_i \Delta u_h^{(i)})^2 / 2. \tag{24}$$

The gradient of discrepancy (Eq. (24)) has the appearance

$$\nabla\varepsilon_k = \frac{\partial\varepsilon}{\partial w_k} = (\Delta u_h^{(0)}, (w_i \Delta u_h^{(i)})) \cdot (\Delta u_h^{(0)}, \Delta u_h^{(k)}). \tag{25}$$

The steepest descent is used for $\{w_i\}$ determination in the form

$$w_k^{n+1} = w_k^n - \tau \nabla \varepsilon_k . \tag{26}$$

The normalization is applied to avoid the exact solution shift past iterations finish

$$\tilde{w}_k^{n+1} = w_k^{n+1} / S^{n+1},\; S^{n+1} = \sum_{i=1}^{N} w_k^{n+1} . \tag{27}$$

The orthogonality criterion ($\varphi$ is the angle between $\Delta u_{\perp}$ and $\Delta u_h^{(0)}$)

$$\varphi = \arccos\left(\frac{(\Delta u_h^{(0)}, w_i \Delta u_h^{(i)})}{\left|\Delta u_h^{(0)}\right| \cdot \left|w_i \Delta u_h^{(i)}\right|}\right) \tag{28}$$

is applied to check the quality of numerical results (convergence of iterations (Eq. (26)). The value of angle $\varphi$ is more informative and lucid if compare with discrepancy $\varepsilon(\vec{w})$ defined by Eq. (24). The completeness of basis $\Delta u_h^{(i)}$ may be implicitly checked by the values of expressions (24) or (28).

The basis $\Delta u_h^{(i)} = A_h^{-1} \delta u_h^{(i)} = u_h^{(i),corr} - u_h^{(i)}$ is computed by the single run of the set of solvers and contains from 2 to 5 numerical solutions. We consider operator $A_h^{-1}$ to be unique for all solutions with the tolerance of the second order.

The estimation of $\Delta u_h^{(i)}$ at first step of the algorithm resembles the defect correction approach [30]. However, the expression $u_{\perp} = w_i u_h^{(i)}$ contains the superposition of true approximation errors. We apply only angles between approximation errors $\Delta u_h^{(i)} = A_h^{-1} \delta u_h^{(i)}$ obtained by the defect correction and we search not for the refined solution (a point in the state space), but for the hypersphere, containing the true solution (its projection on the selected grid).

We use the ensemble of $N+1$ numerical solutions. In order to enhance the reliability we consequently select a solution from the ensemble and define it as the basic solution (centre) $C_h = u_h^{(0)}$. After this step we estimate the radius $R = \left\| u_h^{(0)} - u_{\perp} \right\|$. Finally, we select the maximum radius and the corresponding centre point over all tries as the solution in the Prager&Synge sense.

### 4.3. Nonintrusive methods

It is interesting that the expressions resembling the Eqs. (9), (18) are obtained in [17,18,19] (triangle inequality), [20] (diameter of ensemble), [21] (angle between errors) and have a closely related nature to the Prager&Synge method (enables some estimation of $R > \left\| \tilde{u}_h - u_h^{(i)} \right\|$ for $C_h = u_h^{(i)}$) that we demonstrate below. In general, the approaches by [17-21] may be treated as the variants of the Prager&Synge method for nonorthogonal errors, since the distance between solutions may be used as the majorant of the approximation error under certain additional conditions.

Unfortunately, in the domain of CFD, the approximation errors $\Delta u_h^{(k)}$ are correlated near discontinuities (the errors involve waves with the positive and negative parts) and, so, are not exactly orthogonal. Fortunately, the lack of the rigorous orthogonality may be compensated by some additional information. We consider, herein, the angle between errors, the maximum distance between solutions over the set of independent solutions, and ordering of the errors over their norm (in correspondence with [17-20]).

Let's consider the value of the angle between errors

$$\alpha = \arccos \eta \,, \; \eta = \frac{(\Delta u^{(1)}, \Delta u^{(2)})}{\left\| \Delta u^{(1)} \right\| \cdot \left\| \Delta u^{(2)} \right\|} . \tag{29}$$

Certainly, the inequality (9) is valid at some deviation from errors orthogonality (for the angles between errors greater some critical value $\alpha_*$). However, it is interesting to expand the Prager@Synge approach for the general case, so the analysis for arbitrary angles between errors $\alpha \neq \pi/2$ follows.

### 4.3.1 The error estimation using value of the angle between approximation errors

The value of the angle between approximation errors for $\alpha \neq \pi/2$ determines the relation between the independent and correlated components of the error. We consider two main cases below. For the obtuse angle $\alpha > \pi/2$ (Fig. 3) the elementary geometric considerations demonstrate

$$\left\| u_h^{(1)} - \tilde{u}_h \right\| < \left\| u_h^{(1)} - u_h^{(2)} \right\| . \tag{30}$$

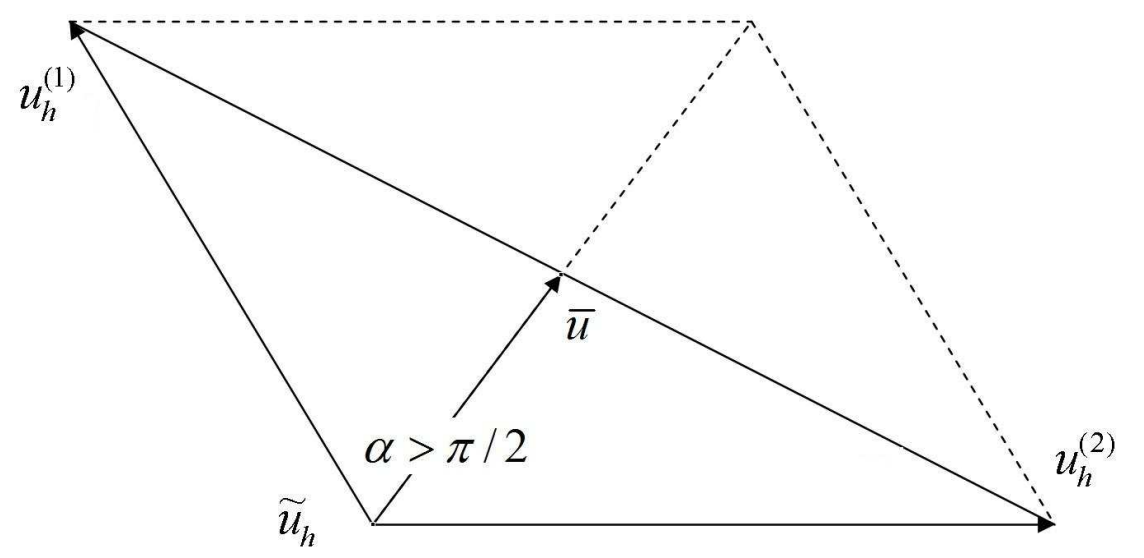

Fig. 3. The geometry of solutions for $\alpha > \pi/2$

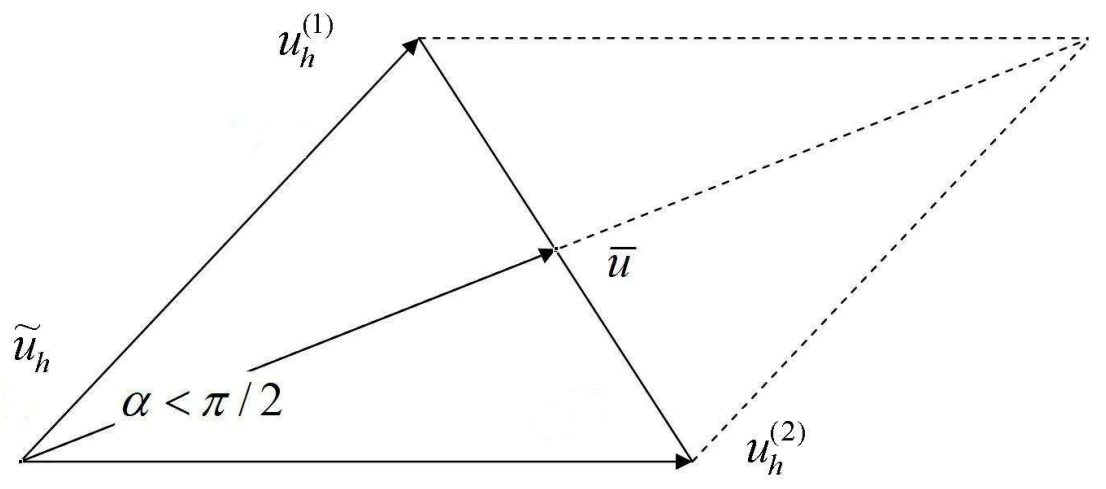

Fig. 4. The geometry of solutions for $\alpha < \pi/2$.

The analysis is more complicated for the acute angle $\alpha < \pi/2$ (Fig. 4). Certainly, the above mentioned inequality (9) is valid for $\alpha$ close to $\pi/2$. Unfortunately, it may break for the small enough value of $\alpha$. The following estimate is valid [21]:

$$\left\| \Delta u^{(k)} \right\| < 1.1 \cdot \frac{\left\| u^{(1)} - u^{(2)} \right\|}{\sin(\alpha/2)} = M, (k = 1,2) . \tag{31}$$

This inequality enables the estimation of the approximation error norm using the value of the angle $\alpha$ between errors and the distance between the numerical solutions. One may see the divergence of this estimate as $\alpha \to 0$ (for parallel errors). So the nonparallelity of errors is crucial for the success of such error estimation. The value of $\alpha$ is unknown but may be evaluated from the analysis of the angles between truncation errors in accordance with [21].

### 4.3.2. The error estimation using triangle inequality

The computable difference of two numerical solutions is equal to the difference between their errors $u_h^{(1)} - u_h^{(2)} = u_h^{(1)} - \tilde{u}_h - u_h^{(2)} + \tilde{u}_h = \Delta\tilde{u}_h^{(1)} - \Delta\tilde{u}_h^{(2)}$ and may be used for the error estimation.

According [17,18] the relation between norms of errors (the detection of most imprecise solution) may be used as an important information at the approximation error estimation. If the norm of difference of two numerical solutions $u_h^{(1)} \in R^M$ and $u_h^{(2)} \in R^M$ $d_{1,2} = \left\| u_h^{(1)} - u_h^{(2)} \right\|$ is known from computations and the relation

$$\left\|\tilde{u}_h - u_h^{(1)}\right\| \geq 2 \cdot \left\|\tilde{u}_h - u_h^{(2)}\right\| \qquad (32)$$

is valid (one of errors is much greater another), then the norm of approximate solution $u_h^{(2)}$ error is less than the norm of the difference of solutions:

$$\left\|\tilde{u}_h - u_h^{(2)}\right\| \leq d_{1,2}\,. \qquad (33)$$

This estimate may be proved by the triangle inequality [18,19,20] for vertexes $u_h^{(1)}, u_h^{(2)}, u_h^{(0)} = \tilde{u}$. For distances (legs of triangle) $d_{0,1} = \left\|\tilde{u}_h - u_h^{(1)}\right\|, d_{0,2} = \left\|\tilde{u}_h - u_h^{(2)}\right\|$, $d_{1,2} = \left\|u_h^{(2)} - u_h^{(1)}\right\|$ the triangle inequality may be presented as $d_{0,1} \leq d_{1,2} + d_{0,2}$ and may be transformed to $d_{0,1} - d_{0,2} \leq d_{12}$. For the expression (32) rearranged into the form $d_{0,2} \leq d_{0,1} - d_{0,2}$, one obtains $d_{0,2} \leq d_{0,1} - d_{0,2} \leq d_{1,2}$ that corresponds to the desired expression (33) $d_{0,2} \leq d_{1,2}$.

The relation (32) may be substantiated at the heuristic level of accuracy by an analysis of the collection of distances between solutions $d_{i,j}$ for the set of three or greater solutions. These distances enable the detection of the nearby and distant solutions. If one of solutions $u_h^{(1)}$ is much less accurate in comparison with others, the set of distances between solutions $d_{i,j}$ splits into a cluster of distances, engendered by the inaccurate solution $u_h^{(1)}$ (great values $d_{1,j}$), and a cluster of distances between more accurate solutions $d_{i,j}(i \neq 1)$. If $u_h^{(1)}$ is sufficiently imprecise ($\min(d_{1,j})$ is twofold (or greater) as $\max(d_{i,j}(i \neq 1))$), the relations (32) and (33) are correct with high probability.

The expression (33) does not assume the error orthogonality. It is based on the *a priori* known relation of error norms for two numerical solutions, and provides the upper estimate only for the more precise solutions (excluding $u_h^{(1)}$).

The selection of the solution maximally distant from other solution is equivalent to the search for the solution which error has the maximum angle with other errors. Thus, the approach by [17,18,19] may be considered as the extension of the Prager&Singe method, based on the analysis of distances between solutions.

### 4.3.3. The error estimation using the width of the ensemble of numerical solutions

Since the exact orthogonality of approximation errors does not occur and above discussed ordering of errors may be not valid, there exists the need for more robust and universal estimates. On the ensemble of $N$ solutions we define the ensemble width (the maximum distance between solutions) [20] as

$$d_{\max} = \max_{i,k}\left\|u_h^{(k)} - u_h^{(i)}\right\| i,k = 1,\ldots,N \qquad (34)$$

and postulate the estimate:

$$\left\|\Delta u_h^{(k)}\right\| \leq d_{\max}\,. \qquad (35)$$

If the approximation errors are orthogonal, the distance between solutions provides the reliable upper estimate of error. If the errors are not orthogonal, the estimate may be violated. The naive idea to circumvent this trouble is to expand the number of numerical solutions $N$ and to select the largest distance between them (this is in some relations similar to the approach described in previous section). As the dimension of the ensemble increases, the smallest distances corresponding the small angles between solutions are neglected and the couple of solutions with the maximum angles (closest to $\pi/2$) between errors are selected. The numerical tests by [20] confirm that the inequality (35) becomes more reliable as the number of the ensemble elements increases. So, the estimation of the width of the ensemble (34) and inequality (35) may be considered as some approximate of the Prager&Synge solution.

## 5. The test problems and numerical algorithms

The above considered methods for the Prager&Synge solution estimation are checked on the numerical solutions for the Euler equations

$$\frac{\partial \rho}{\partial t}+\frac{\partial\left(\rho U^k\right)}{\partial x^k}=0, \qquad (36)$$

$$\frac{\partial\left(\rho U^i\right)}{\partial t}+\frac{\partial\left(\rho U^k U^i+P\delta_{ik}\right)}{\partial x^k}=0, \qquad (37)$$

$$\frac{\partial(\rho E)}{\partial t}+\frac{\partial\left(\rho U^k h_0\right)}{\partial x^k}=0, \qquad (38)$$

describing the flows of the inviscid compressible fluid. The system is two-dimensional, $i,k=1,2$, $U^1=U, U^2=V$ are the velocity components, $h_0=(U^2+V^2)/2+h$, $h=\frac{\gamma}{\gamma-1}\frac{P}{\rho}=\gamma e$, $e=\frac{RT}{\gamma-1}$, $E=\left(e+\frac{1}{2}(U^2+V^2)\right)$ are enthalpy and energies, $P=\rho RT$ is the state equation, $\gamma=C_p/C_v$ is the specific heats relation.

The Edney-I and Edney-VI flow structures [22], presented by the Fig. 5 and Fig. 6, are analyzed due to the existence of analytical solutions used for the estimation of the “exact” discretization error $\Delta\tilde{u}^{(k)}=u^{(k)}-\tilde{u}_h$. This error is obtained by the subtraction of the numerical solution from the projection of the analytical one onto the computational grid.

Fig. 5 provides the spatial density distribution for Edney-I flow pattern ($M=4$, flow deflection angles $\chi_1=20^o$ (upper) and $\chi_2=15^o$ (lower)). Fig. 6 provides the density distribution for Edney-VI ($M=3.5$, consequent angles $\chi_1=15^o$ and $\chi_2=25^o$).

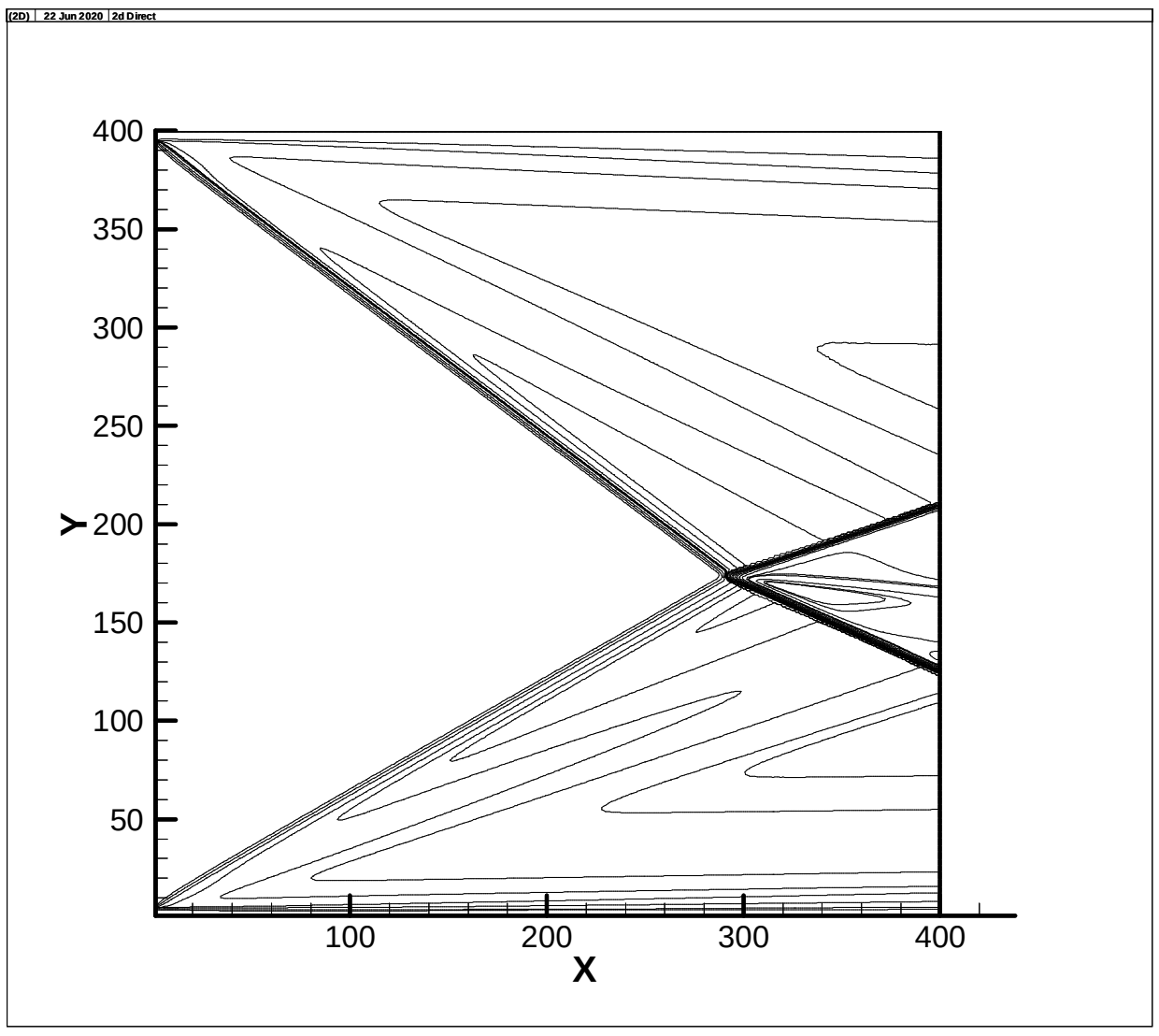

Fig. 5. Density isolines (Edney I)

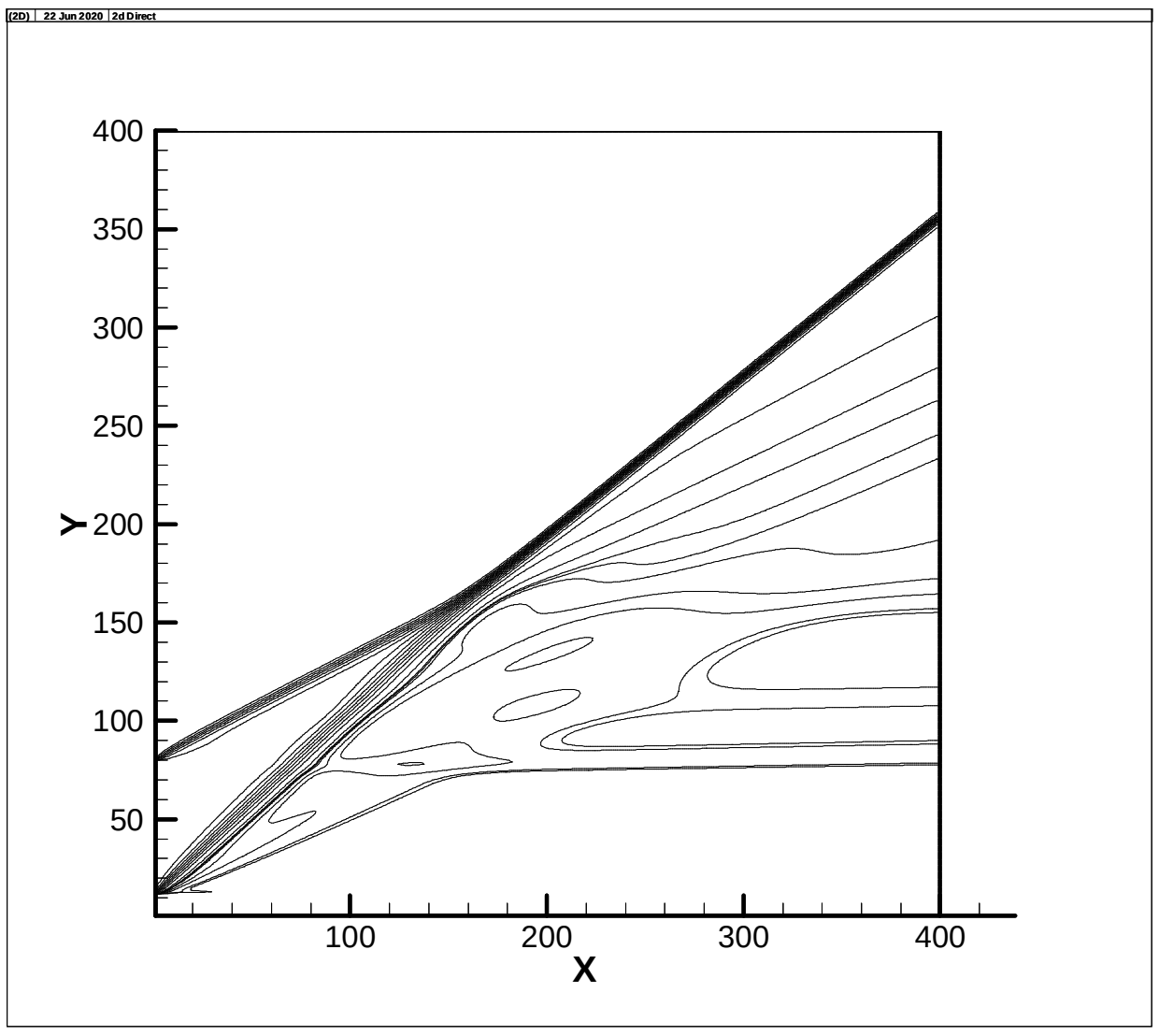

Fig. 6. Density isolines (Edney VI)

We analyze the set of 13 numerical solutions obtained by the numerical methods covering the range of approximation order from one to five. The methods and designations are used that follow.

“S1”: First order algorithm by Courant-Isaacson-Rees [31,32]. It is the most imprecise method, nevertheless, it ensures the success of the “triangle inequality” based estimates.

“S2”: Second order MUSCL [33] based algorithm that uses approximate Riemann solver by [31].

“HLLC”: Also second order MUSCL [33] based algorithm that uses approximate Riemann solver by [46].

These two methods (“S2” and “HLLC”) differ only by the approximation of the Riemann solver, provide very close solutions, and may be used for the worst case estimates. Other methods have essentially different algorithmic structure and potentially are more suitable for *a posteriori* error estimation.

“relax”: Second order relaxation based [38] algorithm [39].

“MC”: Second order MacCormack [40] scheme without the artificial viscosity.

“MC1”: MacCormack scheme with the second order artificial viscosity, (viscosity coefficient $\mu = 0.01$).

“MC2”: MacCormack scheme with the second order artificial viscosity, (viscosity coefficient $\mu = 0.002$).

“MC4”: MacCormack scheme with the fourth order artificial viscosity, (viscosity coefficient $\mu = 0.01$).

“LW”: Second order “two step Lax-Wendroff” method [41,42] with the artificial viscosity of second order ($\mu = 0.01$).

The MacCormack and Lax-Wendroff schemes are rather old and, in the original version, do not provide the monotonicity in the vicinity of discontinuities.

“S3”: Third order algorithm based on the modification of Chakravarthy-Osher method [35,36].

“W3”: Third order algorithm WENO [43].

“S4”: Fourth order algorithm by [37].

“W5”: Fifth order algorithm WENO [44, 45].

## 6. The results of the numerical tests

In numerical tests we obtain the solution in the sense of Prager&Synge. It is determined by the choice of the centre of the hypersphere (some numerical solution) and the radius $R = \|\Delta u_h\|$.

Naturally, these data may be treated as *a posteriori* error estimate. The quality of *a posteriori* error estimation may be described by the effectivity index, which is equal to the ratio of the norm of the error estimate $\|\Delta u_h\|$ to the norm of the true error $\|\Delta \tilde{u}\|_h$

$$I_{eff} = \frac{\|\Delta u_h\|}{\|\Delta \tilde{u}_h\|}. \quad (39)$$

The error estimate should be greater the true error (to ensure the true solution enclosure) and should be not too pessimistic. The condition $I_{eff} \geq 1$ means the successful construction of the Prager&Synge solution. The relation $I_{eff} < 3$ is stated by [7] for the finite-element applications, which imply sufficiently smooth solutions and may be considered as some reference point for CFD applications.

We performed numerical experiments for several cases including different flow structures (Edney-I, Edney-VI), different flow parameters (Mach numbers in the range $3 \div 5$, shock angles in the range $10 \div 30^o$,) and uniform grids of $100\times100$, $200\times200$ and $400\times400$ nodes. In the following subsections several typical results are presented as the illustrations for different methods for computation of the Prager&Singe solution.

The truncation error based estimation (Eqs. (11)-(18)) provides nonrealistic pessimistic results for the efficiency index $I_{eff} \sim 10 \div 20$. In combination with the extreme algorithmic complexity these results make this option non competitive.

The approximation error based estimation (Eqs. (19)-(28)) provides acceptable results described in Section 6.1.

The typical results of the tests for the nonintrusive options (Section 4.3) are presented in Sections 6.2, 6.3, 6.4 for illustration. The results provided by these options are close from the viewpoint of the values of the efficiency index. However, the option based on the ensemble width, does not require the additional information (angles between errors or error relation). By this reason the ensemble widths based method should be preferred in applications.

### 6.1 The estimation of the Prager&Synge solution based on the orthogonalization of approximation errors

Herein, we consider the numerical tests for the algorithm described in the Section 4.2 and based on the auxiliary solution providing the orthogonality of approximation errors. The results for the Edney-VI shock interaction pattern ($M = 3.5$, $\alpha_1 = 15^o$, $\alpha_2 = 25^o$) for the numerical tests on the grids containing $100\times100$ and $400\times400$ nodes are presented as typical.

At first step the solution obtained by [34] (S2) is used as the basic solution $u_h^{(0)}$, which is considered as the centre of the hypersphere containing the true solution (more correctly, the projection of the true solution on considered grid). The additional solutions $u_h^{(1)}$ [31,32], $u_h^{(2)}$ [35], $u_h^{(3)}$ [37], $u_h^{(4)}$ [40], $u_h^{(5)}$ [39] are used in different combinations in order to generate some orthogonal error and the superposition of solutions as an approximation of the auxiliary solution $u_\perp = w_i u_h^{(i)}$. The value $R = \|u_h^{(0)} - u_\perp\|$ is the radius of the hypersphere containing the true solution (the centre is $C_h = u_h^{(0)}$).

At the next steps all other numerical solutions are consequently selected as the new basic solution $u_h^{(0)}$ and the radius $R = \|u_h^{(0)} - u_\perp\|$ is computed.

The minimum value the effectivity index over all basic solutions is above unit that ensures the enclosure of true solution and correctness of the Prager&Synge solution. The maximum value is presented by the Table 1 in the dependence on the number of additional solutions for two different grids.

**Table 1.** The effectivity index in dependence on the number of solutions

| N solutions | 2 | 3 | 4 | 5 |
|---|---|---|---|---|
| $I_{eff}(100\times100)$ | 10.810 | 3.680 | 1.931 | 1.583 |
| $I_{eff}(400\times400)$ | 9.081 | 1.754 | 1.437 | 1.329 |

One may see from the Table 1 that the magnitude of the effectivity index decreases as the number of solutions in use is enhanced. In general, two auxiliary solutions provide too pessimistic estimations of the approximation error. In the range of 3-5 solutions the values of the effectivity index are quite acceptable. Some saturation of the effectivity index is observed at increasing of the number of solutions. The grid resolution rather weakly affects the effectivity index.

### 6.2 The estimation of approximation error using the triangle inequality

We consider the triangle based approach [17,18,19] using four numerical schemes (noted in above Section by "S1-S4") having the different orders of approximation. The results obtained by other algorithms are omitted to avoid the jumble. For Edney-I shock interaction ($M=4$, $\chi_1=10^o$, $\chi_2=15^o$, $100\times100$), the set of distances between solutions splits into clusters $d_{S1,i}$ (distant solutions) and $d_{j,i}$ (close solutions) and $d_{S1,i} > d_{j,i}$. These clusters are engendered by the imprecise solution "S1" (Fig. 7). One may see the successful enclosure of the exact solution in Fig. 8.

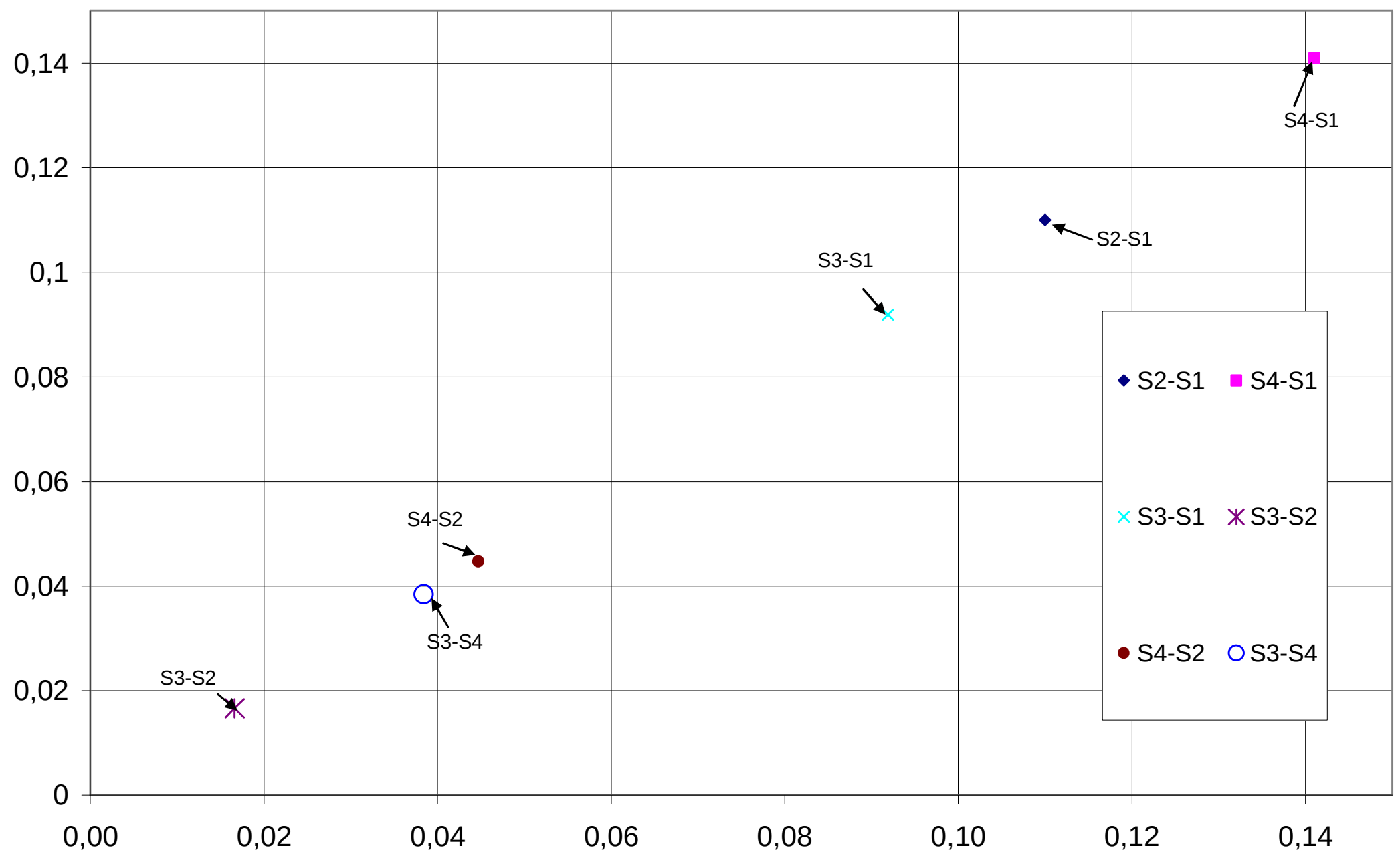

Fig. 7. Clusters of distances between solutions (Edney-I).

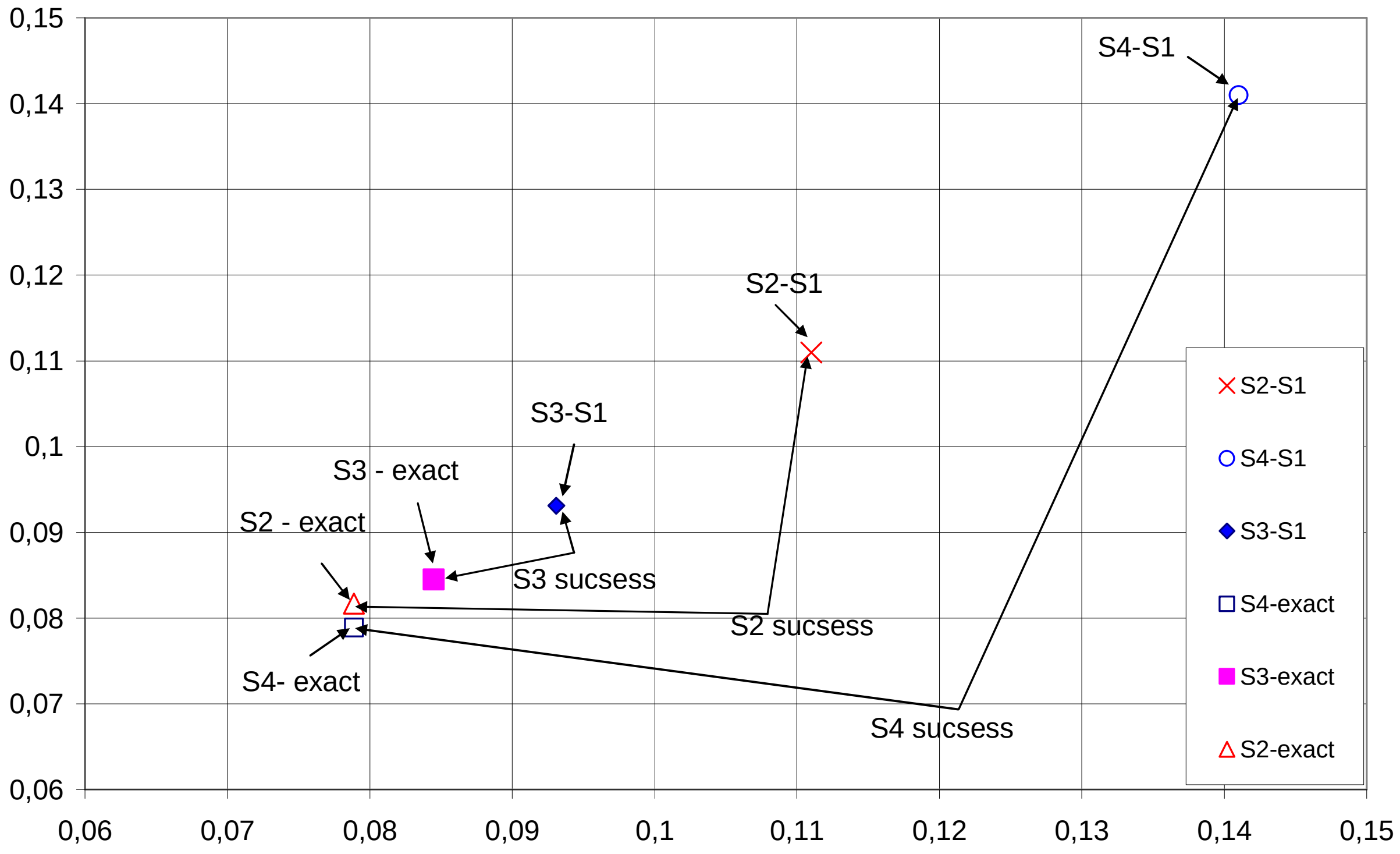

Fig. 8. True error norms vs estimates, obtained using the triangle inequality (Edney-I).

Over all tests, the triangle inequality based approach provides the effectivity index $I_{eff} \sim 0.8 \div 1.5$ if ordering of errors is detected.

### 6.3 The estimation of approximation errors using the angle between them

In the paper [21] the truncation error $\delta u^{(k)}$ was estimated using the postprocessor described by [29]. The angles between truncation errors $\beta_{km}$ are calculated for the ensemble of $N = 13$ solutions obtained using 10 different numerical algorithms [31-46] and several variants of the artificial viscosity. The angles between approximation errors $\alpha_{km}$ are calculated for the same ensemble. This ensemble provides $N \cdot (N-1)/2 = 78$ values of angles for any considered flow patterns. These angles are computed for the ensemble of tests containing two flow patterns and three grids. The averaged (over the ensemble) angle between the truncation errors is in the range $\overline{\beta} = 58 \div 64^o$ for the considered flow patterns and grids. The averaged angle between the approximation errors is in the range $\overline{\alpha} = 30 \div 44^o$. The inequality $\beta_{km} > \alpha_{km}$ holds almost in all tests. The inequality $\alpha_{km} \geq \beta_{km}/3$ is used for the error estimation in accordance with results by [21].

The value of the effectivity index of this estimate is in the range $I_{eff} \sim 0.9 \div 4.5$ (with the mean value $\bar{I}_{eff} \approx 2$) for all considered flow patterns and grids. This approach provides the results of acceptable quality with the minimal computational burden (using only two solutions).

Table 2 presents the effectivity index for the estimates performed by the expression (31) for one of test cases (Edney-I, $M = 4$, $\chi_1 = 20^o$, $\chi_2 = 15^o$, $100 \times 100$).

**Table 2.** The effectivity index of estimate (Eq. (31)) for Edney-I flow mode, grid dimension $100 \times 100$.

| | relax | S2 | HLLC | S3 | S4 | W3 | W5 | MC1 | MC2 | LW | MC | MC4 |
|---|---|---|---|---|---|---|---|---|---|---|---|---|
| S1 | 2.29 | 1.93 | 2.22 | 2.02 | 1.97 | 1.89 | 2.14 | 2.01 | 2.01 | 2.39 | 2.12 | 1.96 |
| relax | 0 | 2.08 | 2.05 | 2.06 | 2.11 | 2.10 | 2.15 | 2.13 | 2.04 | 2.27 | 1.95 | 3.48 |
| S2 | | 0 | 2.37 | 2.20 | 2.22 | 2.31 | 0.92 | 2.22 | 2.09 | 2.25 | 1.99 | 2.09 |
| HLLC | | | 0 | 2.19 | 2.21 | 2.32 | 2.17 | 2.23 | 2.10 | 2.25 | 1.99 | 2.10 |
| S3 | | | | 0 | 2.2 | 2.21 | 2.00 | 2.31 | 2.21 | 2.10 | 2.06 | 2.12 |
| S4 | | | | | 0 | 2.21 | 1.96 | 2.22 | 2.15 | 2.06 | 2.04 | 2.15 |
| W3 | | | | | | 0 | 2.14 | 2.21 | 2.10 | 2.13 | 1.99 | 2.09 |
| W5 | | | | | | | 0 | 1.98 | 1.97 | 2.35 | 2.08 | 1.93 |
| MC1 | | | | | | | | 0 | 2.19 | 2.52 | 2.06 | 2.07 |
| MC2 | | | | | | | | | 0 | 2.12 | 2.25 | 2.13 |
| LW | | | | | | | | | | 0 | 2.1 | 2.05 |
| MC | | | | | | | | | | | 0 | 2.02 |

The inequality $\alpha_{km} \geq \beta_{km}$ is quite natural. Unfortunately, in the inequality $\alpha_{km} \geq C \cdot \beta_{km}$ the constant may be problem dependent.

### 6.4 The estimation of approximation error using the ensemble width

In Table 5 we present the effectivity index value for the error estimations based on the ensemble width (Eq. (33)) ($I_{eff,k} = \max_{i,j}(d_{ij}) / \left\| \Delta u^{(k)} \right\|$) for one of the flow modes (Edney-I, $M = 4$, $\alpha_1 = 20^o$, $\alpha_2 = 15^o$, $100 \times 100$).

**Table 3.** The effectivity index of estimates based on the ensemble width

| | S1 | relax | S2 | HLLC | S3 | S4 | W3 | W5 | MC1 | MC2 | LW | MC | MC4 |
|---|---|---|---|---|---|---|---|---|---|---|---|---|---|
| $\left\| \Delta u^{(k)} \right\|_{L_2}$ | 0.25 | 0.15 | 0.16 | 0.16 | 0.19 | 0.17 | 0.15 | 0.27 | 0.16 | 0.19 | 0.22 | 0.23 | 0.18 |
| $\max_{k,m} d_{km}$ | 0.28 | | | | | | | | | | | | |
| $I_{eff,i}$ | 1.12 | 1.86 | 1.75 | 1.75 | 1.47 | 1.64 | 1.86 | 1.03 | 1.75 | 1.47 | 1.27 | 1.22 | 1.55 |

The effectivity index is in the range $I_{eff} \sim 1.0 \div 2.0$ and is quite acceptable over all considered tests for this (very large) set of solutions. For smaller sets of solutions the quality of results may deteriorate.

The smallest over tests value $I_{eff} = 0.04$ corresponds the couple of AUFS (S1) [34] and HLLC [46] methods, which are intentionally selected as providing the most dependent solutions (both second order accuracy and based on the approximated Riemann problem). In other cases, the results are much more favorable.

### 6.5 The comparative quality of above considered error estimators

The version of the Prager@Synge method based on the adjoint postprocessor and the defect correction (Eq. (18)) overestimates the error and provides the effectivity index in the range $I_{eff} \sim 10 \div 20$ that is too pessimistic.

The version of the Prager@Synge method based on the superposition of solutions (Eq. (19)) provides the effectivity index in the range $I_{eff} \sim 1.0 \div 4.0$.

The use of the angle between truncation errors $\beta_{km}$ and the inequality $\alpha_{km} \geq \beta_{km}/3$ enables the application of the expression (31) with the value of the effectivity index in the range $I_{eff} \sim 0.9 \div 4.5$. This approach provides the results of acceptable quality with the minimal computational burden (using only two solutions). Unfortunately, it requires *a priori* information on the relation of angles between approximation and truncation errors.

A triangle inequality based approach (Eq. (33)) provides the effectivity index $I_{eff} \sim 0.8 \div 1.5$ if ordering of errors over the norm is detected, that is not, unfortunately, the common case.

The effectivity index, based on the ensemble width (Eq. (34)) ($I_{eff,k} = \max_{i,j}(d_{ij})/\left\|\Delta u^{(k)}\right\|$), is in the range $I_{eff} \sim 1.0 \div 2.0$. The value of effectivity index asymptotically improves as the number of solutions increases, however rather slowly. The methods of the close inner structure should be avoided when this approach is used.

For the considered numerical tests the above discussed error norm estimators demonstrate very different laboriousness at close values of effectivity indices. The selection of the suitable method depends on the number and properties of the available solvers.

## 7. Discussion

The classic implementation of Prager@Synge approach is limited by the Poisson or biharmonic equations, which have the specific orthogonality property.

The universal extension of Prager@Synge approach, considered herein, may be applied to the arbitrary system of PDE. As the price, an ensemble of the numerical solutions obtained by independent algorithms is required. In addition, the postprocessors for the truncation error estimation and a defect correction algorithm for the approximation error estimation are necessary. By this reason the universal extension is rather complicated and intrusive. However, it is implemented on single grid (without mesh refinement) and, by this reason, is much more computationally nonexpensive in comparison with the Richardson method.

The orthogonality of truncation errors $\delta u_h^{(k)} \perp \delta u_h^{(i)}$ observed in numerical tests [21] may be the base for the implementation of the first version of the universal approach. Unfortunately, this option is extremely complex form the computational viewpoint and provides too pessimistic results due the instability.

The auxiliary solution $u_\perp = w_i u_h^{(i)}$ having the approximation error orthogonal to the error of some basic solution is another option for the universal extension of Prager@Synge approach. This option is algorithmically simpler and provides results of acceptable quality.

The method, based on the width of the ensemble of numerical solutions, may be considered as an approximate for of the previous option. It is simple from the algorithmic viewpoint (nonintrusive) and computationally nonexpensive.

The Table 4 presents the comparison of classic and universal variants of the Prager@Synge method.

**Table 4**. The classic and universal variants of the Prager@Synge method.

| | Classic | Universal |
|---|---|---|
| 1 | $A\tilde{u} = \rho$, (for $A = B^* B$) | $A\tilde{u} = \rho$ (for any $A$) |
| 2 | $((Bu_h - B\tilde{u}_h),(q - B\tilde{u}_h)) = 0$ | $(u_h - \tilde{u}_h)(u_\perp - \tilde{u}_h) = 0$ |
| 3 | $\|Bu_h - B\tilde{u}_h\| \le \|Bu_h - q\|$ | $\|u_h - \tilde{u}_h\| \le \|u_\perp - u_h\|$ |
| 3.1 | | $u_\perp = A^{-1}(\rho + \delta u_\perp)$, $\delta u_\perp = A_h A_h^* \theta$ Eq. (18) |
| 3.2 | | $u_\perp = w_i u_h^{(i)}$ Eq. (20) |
| | | nonintrusive (ensemble width) |
| 4 | | $\|u_h - \tilde{u}_h\| \le \max\limits_{i,k} \|u^{(k)} - u^{(i)}\| i,k = 1,...,K$ Eq. (34) |

The considered universal extension of the Prager@Synge method may be applied to the solution of arbitrary systems of PDE and can be realized by several algorithms of the different complexity. The availability of the ensemble of numerical solutions obtained by independent algorithms is the only condition that restricts its applicability.

Certainly, the numerical solutions obtained by such algorithms are highly correlated. Fortunately, the truncation and approximation errors for the ensemble of such solutions are far from to be correlated [22,21] and have a significant independent component.

By “independent algorithms” we mean algorithms with the different order of approximation or different structures of solvers (Riemann problem based solvers, finite-differences etc.). Some heuristics for “independent algorithms” may be found in [21]. For the first order PDE the truncation error may be formally expressed as the formal series [28]

$$\delta u^{(k)} = \sum_{m=j_k}^{\infty} C_m^{(k)} h^m \frac{\partial^{m+1}\tilde{u}}{\partial x^{m+1}} \tag{40}$$

where $j_k$ is the approximation order of $k-th$ algorithm. From this standpoint $\delta u^{(k)}$ and $\delta u^{(m)}$ are independent for algorithms of the different order of approximation, since depend on non coinciding set of derivatives. In addition, the transformation $\delta u^{(k)} \to \delta u^{(m)}$ requires the algorithm of the infinite length. This circumstance may be interpreted as the evidence of the algorithmic randomness (“the string is shorter the program”) [47,48]. The algorithmic randomness may be the reason for the truncation error orthogonality observed in [21] due to the geometric properties of the multidimensional spaces [49,50].

## 8. Conclusions

The Prager&Synge method contains the original concept of the numerical solution, which is not based on asymptotics at grid refinement and enables to escape the modern “tyranny" of mesh generation. The highly limited domain of applicability is the main drawback of the Prager&Synge method in the original variant.

The universal (suitable to an arbitrary PDE systems) variant of the Prager&Synge method is proposed that is based on the ensemble of numerical solutions obtained by independent algorithms on the same grid.

The central element of this approach is the solution in the sense of Prager&Synge. Instead of the sequence of solutions, occurring at the grid size diminishing, the Prager&Synge solution deals with the hypersphere containing the projection of the true solution on computational grid. The numerical solution is the centre of this hypersphere. The Prager&Synge solution provides a natural way to *a posteriori* error estimation and natural criterions for mesh refinement termination related with the Cauchy–Bunyakovsky–Schwarz inequality for valuable functionals.

Several versions of realization of the universal Prager&Synge method are discussed. These approaches include intrusive and nonintrusive options.

The intrusive versions are based either on the orthogonality of truncation errors or on the orthogonality of approximation errors.

The version using the truncation errors is extremely complicated from the algorithmic viewpoint. It applies the postprocessors, corresponding a forward (acting on the solution and providing the truncation error estimate) and adjoint (acting on the truncation error) propagators. The estimates of the error norm, obtained using this version, are too pessimistic due to an action of the adjoint postprocessor on the artificial truncation error that may be highly irregular.

The version using the approximation errors applies the superposition of numerical solutions and provides the acceptable values of the efficiency index for the error norm estimation. From the algorithmic viewpoint this version is of the moderate complexity.

The nonintrusive versions are based on the analysis of the ensemble of numerical solutions obtained by independent algorithms. From the algorithmic viewpoint these versions are very simple. The versions, based on the angles between truncation errors, triangle inequality and the width of ensemble of solutions are compared. The approach using the width of ensemble is some simplified, nonintrusive and approximate variant of the intrusive version based on the approximation errors orthogonality.

The Prager&Synge numerical solution is a natural basis for the *a posteriori* error estimation. The intrusive and nonintrusive approaches for the *a posteriori* error estimation are compared for the two dimensional compressible flows with the shock waves. The numerical tests demonstrated the feasibility on the computation of the numerical solution in the sense of Prager&Synge using these approaches.

The algorithms based on the orthogonality of approximation errors (intrusive) and the width of the ensemble of solutions (nonintrusive) are the best options from the efficiency index and algorithmic simplicity viewpoints.